\documentclass[a4paper,11pt]{article}

\usepackage[latin1]{inputenc}
\usepackage[english]{babel}
\usepackage{hyperref}
\usepackage{amsfonts}
\usepackage{amssymb,amsmath,mathrsfs}

\usepackage{bbm}
\usepackage{graphicx,ifthen,calc,color}
\usepackage{mathpazo} % math & rm
\usepackage[scaled]{helvet} % ss
\usepackage{courier} % tt
\normalfont

\usepackage[T1]{fontenc}

\newcommand{\De}{\mathrm{d}} %differential 
\newcommand{\R}{\mathbb{R}} %reals
\newcommand{\Z}{\mathbb Z} %integers
\newcommand{\N}{\mathbb N} %natural numbers
\newcommand{\var}{\mathbf{Var}_N} %Variance in mathrm
\newcommand{\indicator}[1]{\mathbf{1}_{\{#1\}}} %indicator of an event
\newcommand{\cov}{\mathbf{Cov}_N} %Cov in mathrm

\newcommand{\toi}{\to +\infty}
\newcommand{\prob}{\mathbf P}
\newcommand{\av}{\mathbf E} %expectation
\newcommand{\F}{\mathscr F} %sigma -algebra F
\newcommand{\hpts}{\mathcal H_N(\alpha)}

\newcommand{\cnst}{2 \sqrt{2 g}}
\newcommand{\bound}[1]{\partial_2 #1}
\newcommand{\eq}{\begin{equation}}
\newcommand{\eeq}{\end{equation}}
\newcommand{\eqs}{\begin{equation*}}
\newcommand{\eeqs}{\end{equation*}}
\newcommand{\eqa}{\begin{eqnarray}}
\newcommand{\eeqa}{\end{eqnarray}}
\newcommand{\eqas}{\begin{eqnarray*}}
\newcommand{\eeqas}{\end{eqnarray*}}
\newcommand{\aatop}[2]{\genfrac{}{}{0pt}{1}{#1}{#2}}
\def\ci{\perp\!\!\!\perp}

\def\proof{\noindent{\bf Proof.}\ } %command for beginning of the proof
\def\endproof{{\mbox{}\nolinebreak\hfill\rule{2mm}{2mm}\par\medbreak} }
\newcommand{\proofnam}[1]{\noindent{\bf Proof of Theorem \ref{#1}.}\ } %command for beginning of the proof with name

\newtheorem{lemma}{Lemma}[section]

\newtheorem{itrem}[lemma]{Remark}
\newtheorem{itdef}[lemma]{Definition}
\newtheorem{theorem}[lemma]{Theorem}

\numberwithin{equation}{section}

\newcommand{\vr}{\varphi}

\usepackage[all]{xy}
\usepackage{stmaryrd}

\title{High points for the membrane model in the critical dimension}
 \author{Alessandra Cipriani\footnote{Institut f\"ur Mathematik, Universit\"at Z\"urich, Winterthurerstrasse 190, CH-8057, Zurich, Switzerland; \texttt{alessandra.cipriani@math.uzh.ch}. Work supported by the Swiss National Science Foundation, grant 138141, and by the Forschungskredit of the University of Zurich.}}

\begin{document}

\maketitle

\begin{abstract}
In this notice we study the fractal structure of the set of high points for the membrane model in the critical dimension $d=4$. The membrane model is a centered Gaussian field whose covariance is the inverse of the discrete bilaplacian operator on $\Z^4$. We are able to compute the 
Hausdorf\/f dimension of the set of points which are atypically high, and also that of clusters, showing that high points tend not to be evenly spread on the 
lattice. We will see that these results follow closely those obtained by O. Daviaud \cite{Daviaud} for the 2-dimensional discrete Gaussian Free Field.
\end{abstract}

\section{The model}
The field of random interfaces has been widely studied in statistical mechanics. 
These interfaces are described by a family of real-valued random variables indexed by the $d$-dimensional integer lattice, 
which are considered as a height configuration, namely they indicate the height of the interface above a reference hyperplane. 
The probability of a configuration depends on its energy (the Hamiltonian), which defines a measure on the space of such configurations. The most well-known models are the so-called 
\emph{gradient model}, in particular 
the Discrete Gaussian Free Field (DGFF), or harmonic crystal, whose Hamiltonian is a function of the discrete gradient of the heights, and the {\em membrane model}. The study of such interface was firstly undertaken by Sakagawa in \cite{Sakagawa}; we are aware of 
the contributions of Kurt (\cite{Kurt_d4}, \cite{Kurt_d5}) regarding also a phenomenon called \emph{entropic repulsion} in dimension $4$.\\
The Membrane Model is a Gaussian multivariate random variable whose Hamiltonian depends on 
the mean curvature of the interface, in particular favors configurations whose curvature is approximately constant.
It is indeed a lattice-based scalar field $\{\varphi_x\}_{x\in\Z^d}$ where $\varphi_x\in \R$ is viewed as a height variable at the site $x$ 
of the lattice. There are three convenient and equivalent ways in which one can see such a field. Denote by $V_N:=[-N,N]^d\cap\Z^d$ the centered box of side-length $2N+1.$ Then
\begin{enumerate}
\item the membrane model is the random interface model whose distribution is given by
\begin{equation}
 \label{def:field}
\mathbf P_N(\De\varphi)=\frac{1}{Z_N}\exp\left(-\frac{1}{2}\sum_{x\in \Z^d}(\Delta\varphi_x)^2\right)\prod_{x\in V_N}\De\varphi_x
\prod_{x\in \partial_2 V_N}\delta_0(\De\varphi_x),
\end{equation}
where $\Delta$ is the discrete Laplacian, $\partial_2V_N:=\{y\in V_N^c:\,  d(y,V_N)\leq 2\}$ and $Z_N$ is the normalizing constant.
% This means we are imposing 0-boundary conditions. We use the notation
% $$H_N:=-\frac{1}{2}\sum_{x\in \Z^d}(\Delta\varphi_x)^2.$$
\item By re-summation, the law $\mathbf P_N$ of the field is the law of
the centered Gaussian field on $V_N$ with covariance matrix
$$G_N(x,y):=\cov(\varphi_x,\varphi_y)=(\Delta_N^2)^{-1}(x,y).$$
Here, $\Delta_N^2(x,y)=\Delta^2(x,y)1_{\{x,y\in V_N\}}$ is the Bilaplacian with 0-boundary conditions outside $V_N.$ 
\item The model is a centered Gaussian field on $V_N$ whose covariance matrix $G_N$ satisfies, for $x\in V_N,$
$$\left\{\begin{array}{lc}
\Delta^2 G_N(x,y)=\delta_{xy},& y \in V_N \label{eq:uno} \\
G_N(x,y)=0, & y \in \partial_2 V_N.\nonumber
\end{array}\right.$$
\end{enumerate}
For $d\geq 5$ the infinite volume Gibbs measure $P$ exists \cite[Prop. 1.2.3]{Kurt_Thesis} and is the law of the centered Gaussian field
with covariance matrix
$$G(x,y)=\Delta^{-2}(x,y).$$
The membrane model presents several points in common, as well as challenging differences, from the more known DGFF. 
The former lacks some key features of the latter, namely 
\begin{enumerate}
\item the random walk representation for the Green's function. In the harmonic crystal, it is possible to establish the well-known relation involving the covariance matrix $\Gamma_N$:
\begin{equation}\label{eq:harmonic_cov}\Gamma_N(x,\,y)=\mathbb E^x\left(\sum_{n=1}^{\tau_{\partial V_N}-1}\indicator{S_n=y}\right),\end{equation}
where $\mathbb E^x$ is the law of a standard random walk $(S_n)_{n\geq 0}$ started at $x\in \Z^2$ and $\tau_{\partial V_N}$ is the first exit time from $V_N$.
\item Absence of monotonicity, for example the FKG inequality.
\end{enumerate}
It is thus not possible to rely on harmonic analysis to control the field, 
and this renders many problems solved for the harmonic crystal quite intractable. Despite the lack of such tools it is sufficient to establish two crucial properties to study the high points: one is the {\em logarithmic bound} on covariances which are explained in Lemma \ref{lemma2.1}, and the other one is the {\em 2-Markov property}, which can be stated as follows:
\begin{itdef}[$2$-Markov property] Let $A, \,B \subseteq V_N$ and $\mathrm{dist}(A,B)\geq 3$. Then $\left\{ \varphi_x\right\}_{x \in A}$ and $\left\{ \varphi_x\right\}_{x \in B}$ are independent under the conditional law 
$$\mathbf P_N\left(\cdot\, |\,\sigma\left( \left\{ \varphi_x,\,x \notin A \cup B\right\}\right)\right).$$
\end{itdef}
This suggests that the behavior of certain Gaussian fields with respect to exceedences is universal, in the sense that as soon as the model displays a Gibbs-Markov property and covariances decay at the same rate, then the behavior of high points is the same (with some small adjustments to be done according to the dimension). This also opens up the question of whether there are other points in common between log-correlated Gaussian fields, and we believe a more precise answer will be given soon.\\
The starting point is understanding how many ``high'' points viz. points that grow more than the average there are typically. The first step is to find the average height of the field, 
in other words to show that
there exists a constant $c>0$ such that 
$$\mathbf E\left(\max_{x\in V_N}\varphi_x\right)/\log N\stackrel{N\toi}{\longrightarrow} c.$$
\begin{theorem}[{\cite[Theorem 1.2]{Kurt_d4}}]\label{Kurt1.2}
 Let $d=4$, $\ell\in (0,1)$, 
\begin{equation}\label{eq:V_N}V^\ell_N:=\left\{ x\in V_N:\,d(x,V_N^c)\geq \ell N\right\}\end{equation}
 and let $g:=8/\pi^2$. Then
\begin{enumerate}
 \item[(a)] \[\lim_{N \toi}\prob\left(\sup_{x\in V_N}\vr_x\geq \cnst \log N \right)=0. \]
 \item[(b)] If $0<\ell<1/2$, $0<\eta<1$ there exists $C=C(\ell,\,\eta)>0$ such that
  \[ \prob\left(\sup_{x\in V_N^\ell}\vr_x\geq \left(\cnst-\eta\right) \log N \right)\leq \exp\left(-C \log^2 \! N\right).\]
\end{enumerate}
\end{theorem}
Roughly said, the first-order approximation of the maximum is of order $\log N$, which also implies that the field behaves approximately like independent variables.
For us then an $\alpha$\emph{-high point} will be a point whose height is greater than $2\sqrt{2 g}\alpha\log N$. 
The behavior of $\alpha$-high points for the $2$-dimensional DGFF, as shown in \cite{Daviaud}, tells us that such points exhibit a fractal structure. Very similar 
results were obtained by Dembo, Peres, Rosen and Zeitouni in \cite{DPRZ} for the set of late points of the 2-d standard random walk. \\
To begin with, we recall the definition of the discrete fractal dimension:
\begin{itdef}[Discrete fractal dimension, \cite{Barlow}]
 Let $A \subseteq \Z^d$. If the following limit exists, the fractal dimension of $A$ is
\[
 dim(A):=\lim_{N\to \infty}\frac{\log\left|A\cap V_N\right|}{\log N}.
\]
\end{itdef}
The fractal dimension of the high points is given then in
\begin{theorem}[Number of high points]\label{thm1.3}
 Let $\ell\in (0,1)$, and
$$\mathcal H_N(\eta):=\left\{x \in V^\ell_N: \varphi_x \geq 2\sqrt{2 g}\eta \log N \right\}$$
be the set of $\eta$-high points.
\begin{enumerate}
 \item[(a)] For $0<\eta<1$ we obtain the following limit in probability:
\eq\nonumber
 \lim_{N \toi} \frac{\log |\mathcal H_N(\eta)|}{\log N}=4(1-\eta^2).
\eeq
\item[(b)] For all $\delta>0$ there exists a constant $C>0$ such that for $N$ large
\[
 \mathbf P_N\left(\left\{|\mathcal H_N(\eta)|\leq N^{4(1-\eta^2)-\delta} \right\}\right)\leq \exp(-C\log^2 \! N).
\]
\end{enumerate}
\end{theorem}
We can push further the comparison between the DGFF and the Membrane Model at their respective critical dimensions, and one can find an 
interesting similarity in the behavior of the points. \cite{Daviaud} for example also showed that high points appear in clusters; this is what 
occurs in the membrane model, as the following two theorems show:
\begin{theorem}[Cluster of high points 1]\label{thm1.4}
 Let 
$$D(x,\rho):=\left\{y \in V_N:\,|y-x|\leq \rho \right\}.$$
For $0<\alpha<\beta<1$ and $\delta>0$
\eq \label{eq:1.7}
\lim_{N \toi}\max_{x \in V_N^\ell}\mathbf P_N\left(\left|\frac{\hpts \cap D(x,N^\beta)}{\log N} -4\beta(1-(\alpha/\beta)^2)\right|>\delta \right)=0.
\eeq
\end{theorem}
\begin{theorem}[Cluster of high points 2]\label{thm:cluster2}
 For $0<\alpha<1$, $0<\beta<1$ and $\delta>0$ we have
\[ 
 \lim_{N \toi} \max_{x \in V_N^\ell}\mathbf{P}\left(\left|\frac{\log |\mathcal H_N(\alpha) \cap D(x,N^\beta)|}{\log N} -4 \beta(1-\alpha^2)\right|>\delta\left| 
x \in \mathcal H_N(\alpha)\right.\right)=0.
\]
\end{theorem}
It is also possible to evaluate the average number of pairs of high points as in the following theorem:
\begin{theorem}[Pairs of high points]\label{thm1.6}
Let $0<\alpha<1$, $0<\beta<1$ and let
\[F_{h,\beta}(\gamma):=\gamma^2(1-\beta)+\frac{h(1-\gamma(1-\beta))^2}{\beta}
\]
\[\Gamma_{\alpha,\beta}:=\left\{\gamma\geq 0:\,4-4\beta-4\alpha^2 F_{0,\beta}(\gamma)\geq 0 \right\}=\left\{\gamma\geq 0:\,(1-\alpha^2\gamma^2)\geq 0\right\},
\]
\[
\rho(\alpha,\beta):=4+4\beta -4\alpha^2\inf_{\gamma\in \Gamma_{\alpha,\beta}}F_{2,\beta}(\gamma)>0.
\]
Note that $\Gamma_{\alpha,\beta}=[0,\,1/\alpha]$ is independent of $\beta$.
Then the following limit in probability holds:
\[\lim_{N\toi}\frac{\log\left|\left\{(x,y)\in \hpts:\,|x-y|\leq N^\beta \right\} \right|}{\log N}=\rho(\alpha,\beta) .
\]
\end{theorem}
Finally we can also show what the maximum width of a spike of given length is:
\begin{theorem}[The biggest high square]\label{thm:highsquare}
Let $-1<\eta<1$, $D_N(\eta)$ the side length of the biggest sub-box for which all height variables are uniformly greater than $2 \sqrt{2 g} \eta \log N$, i. e.
\[ D_N(\eta):=\sup\left\{a \in \N: \,\exists x \in V_N^\ell: \,\min_{y \in B(x,a)}\varphi_y\geq 2 \sqrt{2 g} \eta \log N\right\}.
\]
Then the following limit in probability holds:
\[ \lim_{N \toi} \frac{\log D_N(\eta)}{\log N}=\frac{1-\eta}{2} .
\]
\end{theorem}
% \begin{itrem}
%  The following proofs follow closely those of \cite{Daviaud}. This is due to the fact that the basic ``ingredients'' on which they rely are the 
% logarithmic behavior of covariances and the multiscale decomposition of the field, features that hold both for the DGFF in $\Z^2$ and the membrane model 
% in $\Z^4$ respectively.
% \end{itrem}
The paper is organized as follows: in Section 2 we will prove some preliminary results that will be used for the proofs of the main theorems, to which Section 3 is going to be devoted. 
\section{Preliminary Lemmas and results}
\subsection*{Notation}
$D(x,a)$ (resp. $D(x,a]$) denotes the open (resp. closed) Euclidean ball of center $x$ and radius $a$, while $B(x,a)$ is a box centered at $x$ of side length $a$. For the rest of this notice, recall the definition \eqref{eq:V_N} and we let once and for all $\ell \in (0,1/2)$. Let $x_0 \in V_N$ and
$$
M_\alpha:=\left\{x_0+i(N^\alpha+4):\,i\in N^4 \,\mathrm{and}\,x_0+i(N^\alpha+2)\subset V_N\right\}.
$$
We denote by
$x_B$ the center of a (sub)box $B$ and as $\Pi_\alpha$ the union of sub-boxes of side-length $N^\alpha$ 
(without discretization issues) and midpoint in $M_\alpha$. $\F_\alpha$ will be the sigma-algebra generated by $\left\{\vr_x\right\}$ for 
$x \in \bigcup_{B\in \Pi_\alpha}\partial_2 B$. 
Practically we denote with $\Pi_\alpha$ a set of disjoint boxes separated by layers of thickness 2, which thanks to the 2-Markov property will enable us to perform a decomposition procedure on these sets.\\
Furthermore $\varphi_B:=\mathbf E\left(\varphi_{x_B}|{\F_{\partial_2 B}}\right)$ and $\mathbf{Var}_B(\vr_x):=\var\left(\vr_x|{\F_{\partial_2 B}}\right)$.
\subsection{Lemmas}
\subsubsection{The function $\overline G_N(\cdot,\,\cdot)$}
%%%%%%%%%%%%%%%%%%%%%%%%%
%%%Function \overline G
%%%%%%%%%%%%%%%%%%%%%%%%%
In order to prove some of the next results we will introduce the convolution of the harmonic Green's function, which will prove to be a key tool to obtain the crucial estimates on the covariances of our model. 
Let $A$ be an arbitrary subset of $\Z^4$, and for  $x\in A$ let $\Gamma_A(x,\cdot)$ be the solution of the discrete boundary value problem
$$\left\{\begin{array}{lc}
\Delta \Gamma_A(x,y)=\delta_{xy},& y \in A  \\
\Gamma_A(x,y)=0, & y \in \partial A.\end{array}\right.$$
Note that $\Gamma_N$ as in \eqref{eq:harmonic_cov} is the unique solution to the above problem for $A:=V_N$. The convolution of $\Gamma_N$ is 
$$\overline G_N(x,y):=\sum_{z\in V_N}\Gamma_N(x,z)\Gamma_N(z,y),\quad x,\,y\in V_N. $$
\cite{Kurt_d4} contains several bounds and properties of such a function, and we would like here to recall those that we are going to use in the sequel: for all $x,\,y \in V_N$
\begin{itemize}
 \item symmetry: $\overline G_N(x,y)=\overline G_N(y,x)$, 
%\item radial function: $\overline G_N(x,y)=\overline G_N(0,x-y)$,
\item \cite[Lemma 2.2]{Kurt_d4} if $\ell\in (0,1/2)$ there exist $c_1=c_1(\ell)>0$, $c_2>0$ such that \eq\label{Noemi_lemma_2.2}g \log N+c_1\leq \overline G_N(x,y)\leq g \log N+c_2  \eeq
\end{itemize}
With this in mind it is now easier for us to show how to bound the variances and covariances of our field.
%%%%%%%%%%%%%%%%%%%%%%%%%
%%%Bounds on covariances
%%%%%%%%%%%%%%%%%%%%%%%%%
\begin{lemma}[Bounds on the variances]\label{lemma2.1}
Let $d=4$ and $0<\delta<1$. Then
\begin{itemize}
\item there exists $C>0$ such that \eq\label{eq:var}\sup_{x \in V_N}\var(\varphi_x)\leq g \log N+C. \eeq
\item There exists $C(\ell)>0$ such that \eq\label{eq:absvar}\sup_{x \in V_N^\ell}|\var(\varphi_x)- g \log N|\leq C(\ell).\eeq
\item There exist $C>0$ and $C(\ell)>0$ such that
\eqa\label{eq:cov} &&\sup_{\aatop{x, y \in V_N^\ell}{ x \neq y}} \cov(\varphi_x,\vr_y)- g (\log N-\log|x-y|)\leq C.
\eeqa
\eqa\label{eq:abscov} &&\sup_{\aatop{x, y \in V_N^\ell}{ x \neq y}} \left|\cov(\varphi_x,\vr_y)- g (\log N-\log|x-y|)\right|\leq C(\ell).
\eeqa
\end{itemize}
\end{lemma}
\proof For the variances see \cite[Proposition 1.1]{Kurt_d4}. For the covariances, remember that in \cite[Corollary 2.9]{Kurt_d4} that for all $d\geq 4$ and for all $x \in V_N^\ell$
\eq\label{eq:diff_cov} \sup_{y \in V_N^\ell}|G_N(x,y)-\overline G_N(x,y)|\leq c=c(\ell)<+\infty.
\eeq
It is therefore suf\/ficient to show that \eqref{eq:cov} and \eqref{eq:abscov} hold for $\overline G(\cdot, \cdot)$. But we have from \cite[Lemma 2.10]{Kurt_d4}, that there exists a constant $K$ 
such that in $d=4$ for $x\neq y$ and all $\alpha \in (0,2)$
\[ \overline G_N(x,x)-\overline G_N(x,y)=g\log|y-x|+K+o\left(|y-x|^{-\alpha}\right).\]
Hence
\eqas&&
G_N(x,\,y)\leq \overline G_N(x,\,y)+c=  \overline G_N(x,\,x)-g\log|y-x|+K'\stackrel{\eqref{Noemi_lemma_2.2}}{\leq}\\
&&\leq g\log N -g\log|y-x|+K'.
\eeqas
The other bound follows similarly by considering \eqref{eq:abscov}.\endproof
%%%%%%%%%%%%%%%%%%%%%%%%
%%% decomposition of the field
%%%%%%%%%%%%%%%%%%%%%%%%%
Next we give a decomposition of the field which is similar to the one existing for the DGFF (see for example \cite[Section 2.1]{ASS}). With this in mind, we can prove that conditioning on the values of the field assumed on the double boundary of a subset of $V_N\subseteq \Z^4$ (in fact of any $\Z^d$) the resulting field is again the membrane model restricted to the interior of the smaller domain.
\begin{lemma}\label{lemma:Georgii}
Let $B \subseteq V_N$. Let $\F:=\sigma(\vr_z,\,z \in V_N\setminus B)$. Then
\[
 \{\vr_x\}_{x \in B} \stackrel{d}{=} \left\{\mathbf E_N\left[\vr_x|\F\right]+ \psi_x\right\}_{x \in B}
\]
where ``$\stackrel{d}{=}$'' indicates equality in distribution, in particular under $\mathbf P_N(\cdot)$
\begin{enumerate}
 \item[(a)] $\psi_x \ci \F$;
 \item[(b)] $\{\psi_x\}_{x \in B}$ is distributed as the membrane model with 0-boundary conditions on B.
\end{enumerate}
\end{lemma}
\proof Set $\psi_x:=\vr_x-\mathbf E\left[\vr_x|\F\right]$ for all $x \in B$. We have to show that the above results hold.
\begin{enumerate}
 \item[(a)] It is clear from the definition.
\item[(b)]  
Being $\mathbf P_N$ a Gibbs measure, it satisfies the DLR equation: for all $A \subseteq V_N$, $\F_{A^c}:=\sigma\left(\vr_z,\,z \in A^c\right)$,
\eq \label{eq:DLR}
\mathbf P_N(\cdot\,|\,\F_{A^c})(\eta)=P_{A,\eta}(\cdot)\quad \mathbf P_N(\De \eta)- a. \, s.
\eeq
with
\[
 P_{A,\eta}(\De\varphi)=\frac{1}{Z_A}\exp\left(-\frac{1}{2}\sum_{x\in \Z^d}(\Delta\varphi_x)^2\right)\prod_{x\in A}\De\varphi_x
\prod_{x\in V_N\setminus A}\delta_{\eta_x}(\De\varphi_x).
\]
In other words, $P_{A,\eta}$ is a Gaussian distribution with covariance matrix $\left(\Delta_A^2 \right)^{-1}$. Since $\cov(\cdot,\cdot|\F_{A^c})$ 
%is independent of the boundary values $\psi$\footnote{A general fact which for instance is proved in \cite{AdlerTaylor} is that 
%given a Gaussian vector $(\vec{x}_1,\, \vec{x} _2)\sim \mathcal N\left( \vec \mu:=(\vec{\mu}_1,\vec{\mu}_2),\Sigma:=\left(\begin{array}{cc}\Sigma_{11} & \Sigma_{12}\\ \Sigma_{21} & \Sigma_{22}\end{array}\right)\right)$ 
%the distribution of $ x_i$ conditional on $ x_j$ is again Gaussian with mean resp. covariances given by
%\begin{eqnarray}
%\mu_{i|j}=\mu_i+\Sigma_{ij}\Sigma^{-1}_{jj}(x_j-\mu_j)^T \label{eq:conditionalmean}\\
%\Sigma_{i|j}=\Sigma_{ii}-\Sigma_{ij}\Sigma_{jj}^{-1}\Sigma_{ji}. \label{eq:conditionalcov}
%\end{eqnarray}
%Hence the covariance matrix does not depend on the value of $x_j$.} 
we find out that 
it equals $G_A$. In our case this means that $\cov(\cdot|\F)$ is deterministic and equal to $G_B$. So
%\eqas &&
$$\cov(\psi_x,\,\psi_y)=\cov\left(\psi_x,\,\psi_y|\F\right) =\cov\left(\vr_x,\,\vr_y|\F\right)=G_B(x,y)$$
%\eeqas
\end{enumerate}\endproof
%We now have from \cite[Prop. 13.13]{Georgii} that 
%\eq \label{eq:Georgii}
% \mathbf E\left[\vr_x|\F\right]=-\sum_{h \in B}G_B(x,h)\sum_{k \in \Z^d\setminus B}\Delta^2(h,k)\vr_k.
% \eeq
%Hence
%\eqas &&\cov\left(\mathbf E\left[\vr_x|\F\right],\vr_y|\F\right) =\\
%&& =-\cov\left(\sum_{h \in B}G_B(x,h)\sum_{k \in \Z^d\setminus B}\Delta(h,k)\vr_k,\vr_y|\F\right)=\\
%&& = -\sum_{h \in B}G_B(x,h)\sum_{k \in \bound{B}}\Delta^2(h,k) G_B(k,y)=0.
%\eeqas
%since $G_B(k,y)=0$ for $k \in \bound{B}$.

\begin{itrem}
 This result gives us a decomposition of the membrane model in all dimensions.
\end{itrem}

\begin{lemma}\label{lemma8.2}
Let $0<\alpha<1$ and $0<\beta<1$, $\delta>0$ and we define
\[ S=S(\epsilon):=\left\{(x,y)\in V_N^\ell:\,N^{\beta(1-\epsilon)}\leq |x-y|\leq N^\beta \right\}.
\]
Then there exist $C$, $\epsilon_0>0$ (which can be chosen uniformly on $(\alpha,\beta)$ on compact sets of $(0,1)^4$) and $\gamma_\star:=2(2-\beta)^{-1}$ such that for all $\epsilon\leq \epsilon_0$ and all $N$
\[ \max_{(x,y)\in S}\prob\left(x, y \in \hpts \right)\leq CN^{-4\alpha^2 F_{2,\beta}(\gamma_\star)+\delta}.
\]
\end{lemma}
\proof Let $Z:=\vr_x+\vr_y$ and we see that 
$$\left\{x, y \in \hpts\right\}\subseteq \left\{Z\geq 4 \sqrt{2 g} \alpha\log N \right\}.$$
We obtain also 
from \eqref{eq:cov} that
\[\cov(\vr_x,\vr_y)\leq g \log N-g\beta(1-\epsilon)\log N+O(1).
\]
Thus by \eqref{eq:diff_cov} and \eqref{eq:var} 
$$\var(Z)\leq (2g(2-\beta)+O(\epsilon) + O(1/\log N))\log N .$$
Since $F_{2,\beta}(\gamma_\star)=\gamma_\star$, using \eqref{eq:2.6}
\eqas && \prob(Z\geq 4 \sqrt{2 g}\alpha \log N)  \leq \\
& \leq & \exp\left(-\frac{16(\sqrt{2 g})^2\alpha^2\log^2 \! N}{2((2g(2-\beta)+O(\epsilon) + O(1/\log N))\log N} \right)\leq\\
& \leq & \exp\left(-4\alpha^2\gamma^\ast (1+O(\epsilon)+O(1/\log N))\log N\right)\leq \\
& \leq & CN^{-4\alpha^2 F_{2,\beta}(\gamma_\star)+O(\epsilon)}.
\eeqas
\endproof

\begin{lemma}\label{lemma8.1}
 Let $B:=B(x,4N^\beta)$, $\epsilon >0$, $b^\pm(\alpha,\beta,\epsilon,N)=\cnst(\alpha(1-\beta)\pm\epsilon)\log N $, $I(\alpha,\beta,\epsilon,N):=[b^-(\alpha,\beta,\epsilon,N),
b^+(\alpha,\beta,\epsilon,N)]$. Then
\[
 \max_{x \in V^\ell_N}\prob(\varphi_B\notin I(\alpha,\beta,\epsilon,N))|\varphi_x \geq \cnst \alpha \log N)\stackrel{N \toi}{\longrightarrow} 0.
\]
\end{lemma}
\proof We shorten $I$, $b^+$ and $b^-$ for the above quantities. We recall here two useful facts about normal random variables (whose short proof is postponed to the appendix). If $X\sim \mathcal N(0,1)$ then
\begin{equation} \label{eq:2.6}
 P(|X|\geq a)\leq \exp(-a^2/2),\quad \forall\,a \geq 0,
\end{equation}
\eq\label{eq:2.15}
P\left(|X|\geq a\right)\geq \frac{\exp(-a^2/2)}{\sqrt{2 \pi}a},\quad \forall\,a \geq 1.
\eeq
For $\eta>0$ we obtain with \eqref{eq:2.6} and \eqref{eq:2.15}
\[\prob(\varphi_x\geq \cnst \alpha(1+\eta)\log N|\varphi_x\geq \cnst \alpha\log N) \rightarrow 0.\]
as $N \toi$.
This yields
\eqas 
&& \prob(\varphi_B\notin I|\varphi_x\geq \cnst \alpha\log N)=o(1)+\\
&& +\prob(\varphi_B\notin I,\varphi_x\leq \cnst \alpha(1+\eta)\log N|\varphi_x\geq \cnst \alpha\log N)\leq \\
&& \leq o(1)+\prob(\varphi_B\notin I|\varphi_x\in (1,1+\eta)\cnst \alpha\log N).
\eeqas
Now we write $\varphi_x=\varphi_x-\varphi_B+\varphi_B$ and observe that $\varphi_B \ci \varphi_x-\varphi_B$. Therefore $\cov(\varphi_x,\varphi_B)=\var(\varphi_B)$ and so there exists $Z\sim \mathcal N(0,\sigma^2_Z)$, $\sigma^2_Z>0$, for which
\[ \varphi_B=\frac{\var(\varphi_B)}{\var(\varphi_x)}\varphi_x+Z,  \quad Z\ci \varphi_x.
\]
If $x$ is the center of $B \subseteq C$ we can decompose the variances as $\mathbf{Var}_C(\varphi_x)=\mathbf{Var}_C(\varphi_B)+\mathbf{Var}_B(\varphi_x)$, and with this
\[
\frac{\var(\varphi_B)}{\var(\varphi_x)}=(1-\beta)+O\left(\frac{1}{\log N}\right).
\]
It must then be that $\var(Z)=O(\log N)$. Consequently
\eqas &&
\prob(\varphi_B \geq b^+|\varphi_x\in (1,1+\eta)\cnst \alpha\log N)\leq\\
&&\leq \prob \left(Z+\left((1-\beta)+O\left(\frac{1}{\log N}\right)\right)(1+\eta)\cnst \alpha\log N \geq b^+\right)\to 0
\eeqas
for $\eta<\epsilon/(\alpha(1-\beta))$. Similarly
\eqas &&
\prob(\varphi_B \leq b^-|\varphi_x\in (1,1+\eta)\cnst \alpha\log N)\leq\\
&&\leq \prob \left(Z+\left((1-\beta)+O\left(\frac{1}{\log N}\right)\right)\cnst \alpha\log N \leq b^-\right)\to 0.
\eeqas \endproof
%%%%%%% Lemma inutile qui
%%%%%%%%%
% We now mention an estimate on the number of centers of the box structure
% \begin{lemma}\label{lemma8}
% Let $1/2 < \delta<1$, $1/2 < \alpha<1$. Then there exist $\kappa=\kappa(\alpha)$, $a=a(\alpha,\delta)$ such that
% \[
% P(|\left\{B \in \Gamma_\alpha:\varphi_B\geq 0 \right\}|\leq N^\kappa)\leq \exp(-a \log^2 \! N).
% \]
% \end{lemma}
% \proof It is the content of \cite[Lemma 3.2]{Kurt_d4}. \endproof
% \begin{itrem}
% Noemi avoids \cite[Lemma 13]{BDG} in favour of a proof with Sobolev norms B.4.1.
% \end{itrem}

\begin{lemma}\label{lemma8.3}
We keep the notation of Lemma \ref{lemma8.2}. Let $0<\alpha<\beta<1$ and $\delta>0$. For $(x,y)\in S$ define $T(x,y)$ as the set of sub-boxes of side length $2 N^\beta$ such that the centered subbox of side length $N^\beta$ 
contains $x,y$. Then we can find $C,\epsilon_0>0$ such that for $\epsilon\leq \epsilon_0$ and all $N$
\eqas && \max_{\aatop{x,y \in S}{B\in T(x,y)}} \prob\left(\left\{x,y \in \hpts\right\}\cap \left\{ \vr_B \leq\cnst \alpha \gamma(1-\beta)\log N\right\} \right)\\
&&\leq C N^{-4\alpha^2 F_{2,\beta}(\min\left\{ \gamma,\gamma_\star\right\})+\delta}.\eeqas
$\epsilon_0$ can be chosen uniformly on $(\alpha,\beta)$ on compact sets of $(0,1)^4$.
\end{lemma}
\proof
Define
\[E:=\left\{x,y \in \hpts\right\}\cap \left\{ \vr_B \leq\cnst \alpha \gamma(1-\beta)\log N\right\} .\]
We distinguish two cases:
\begin{description}
 \item[$\gamma\geq \gamma_\star$.] We have $\prob(E)\leq \prob(\left\{x,y \in \hpts\right\})$: the claim follows from Lemma \ref{lemma8.2} because $\min\left\{\gamma,\gamma_\star\right\}=\gamma_\star$.
  \item[$\gamma< \gamma_\star$.] It follows from the definition of $\gamma_\star$ that $\gamma< \gamma_\star$ implies $\gamma<2(2-\beta)^{-1}$. For this reason set $a:=1-\gamma(1-\beta)>0$ and $b:=\gamma(2-\beta)-2<0$. 
Letting $Z:=a(\vr_x+\vr_y)+b\vr_B$
\[E\subseteq \left\{ Z\geq (2a+b\gamma(1-\beta))\alpha\cnst \log N\right\}. \]
Furthermore we have the usual decomposition
\eqa&& \var(Z)=a^2\var(\vr_x) +a^2\var(\vr_y)+b^2\var(\vr_B)+\nonumber\\
&&+ 2ab\cov(\vr_x,\vr_B)+2ab\cov(\vr_y,\vr_B)+\nonumber\\
&&+2a^2\cov(\vr_x,\vr_y).\label{eq:8.6}\eeqa
By Lemma \ref{lemma2.1}
\[\var(\vr_B)=\var(\vr_{x_B})-\mathbf{Var}(\vr_{x_B}|{\F_{\partial_2 B}})\leq g(1-\beta)\log N+O(1).\]
and
\eqas&&\cov(\vr_x,\vr_B)=\av(\av(\vr_x|\F_{\partial_2 B})\av(\vr_{x_B}|\F_{\partial_2 B}))=\\
&&=\cov(\vr_x,\vr_{x_B})-\mathbf{Cov}(\vr_x,\vr_{x_B}|{\F_{\partial_2 B}})\geq\\
&&\geq g(\log N-\log|x-x_B|)-g(\beta\log N-\log|x-x_B|)+O(1)=\\
&&=g(1-\beta)\log N+O(1).
\eeqas
Analogously
\[\cov(\vr_y,\vr_B)\geq g(1-\beta)\log N+O(1).\]
Define the auxiliary function $f(a,b,\beta):=2a^2(2-\beta)+b^2(1-\beta)+4ab(1-\beta)$. We use these bounds in \eqref{eq:8.6} to obtain
\[\var(Z)\leq (f(a,b,\beta)+O(\epsilon)+O(1/\log N))g\log N.\]
By the equality $2a+b=\gamma\beta $
\[4a^2+b^2+4 ab=(2a+b)^2=\gamma^2 \beta^2.\]
Then
\eqas f(a,b,\beta) &=& (2a+b)^2-\beta(2a^2+b^2+4ab)=\\
&=&(4a^2+b^2+4 ab)(1-\beta)+2\beta a^2= \\
&=& (\gamma \beta)^2(1-\beta)+2\beta a^2=\\
&=& \beta(\beta \gamma^2(1-\beta)+2 a^2)=\\
&=& \beta( (2a+b)(1-a)+2 a^2)=\\
&=& \beta(2a+b-ab).
 \eeqas
Hence
\eq\label{eq:8.8}
\var(Z)\leq (\beta(2a+b-ab)+O(\epsilon)+O(1/\log N))g \log N.
\eeq
Since $2a+b-ab=2a+b\gamma(1-\beta)$ \eqref{eq:8.6} and \eqref{eq:8.8} yield
\[\prob(E)\leq C\exp\left(-\left(\frac{4\alpha^2(2a+b-ab)}{\beta}+O(\epsilon)\right)\log N \right). \]
Finally notice that
\eqas && \beta F_{2,\beta}(\gamma)=\beta\gamma^2(1-\beta)+2(1-\gamma(1-\beta))^2=\beta\gamma^2(1-\beta)+2a^2=\\
&&=(2a+b)(1-a)+2 a^2=2a+b-ab.
\eeqas
This allows us to conclude the proof.\endproof
Finally we would like to recall 
\end{description}
\begin{lemma}[{\cite[Lemma 2.11]{Kurt_d4}}]\label{lemma12}
 Let $0<n<N$, $A_N\subseteq \Z^4$ be a box of side-length $N$, $A_n\subseteq A_N$ a box of side-length $n$. %with the 
%same center $x_B \in \Z^4$ as $A_N$ 
Let $0<\epsilon<1/2$. There exists $C>0$ such that for all $x
\in A_n$ with $|x-x_B|<\epsilon n$
\[
 \mathbf{Var_N}\left(\mathbf{E}\left( \vr_x\,|\,\F_{\bound{A_n}}\right)-\mathbf{E}_N\left( \vr_{x_B}\,|\,\F_{\bound{A_n}}\right)\,|\,\F_{\bound{A_N}}\right)\leq C \epsilon
\]
\end{lemma}

\section{Five theorems}
% \begin{description}
%  \item[Upper bound.] We have
% \begin{eqnarray*}
% && \prob\left(|\hptseta|\geq N^{4(1-\eta^2)+\delta}\right)\stackrel{Chebyshev\,inq.}{\leq} N^{-4(1-\eta^2)-\delta}\mathbf E(|\hptseta|)\leq\\
% && \leq N^{-4(1-\eta^2)-\delta}N^4 \max_{x \in V_N}\prob(\varphi_x\geq \cnst \eta \log N)\leq\\
% && \leq N^{-4(1-\eta^2)-\delta}N^4 \exp{\left[-\frac{8 g \eta^2 \log^2 \! N}{2 g \log N +c_3} \right]}\stackrel{N \toi}{\longrightarrow} 0
% \end{eqnarray*}
% \item[Lower bound.] \cite[Theorem 3.1]{Kurt_d4}.
% \end{description}
% \endproof
\proofnam{thm1.3} The core of the proof is the lower bound (b) which was already proved by \cite[Theorem 1.3]{Kurt_d4} and is based on the hierarchical 
decomposition of the membrane model, similar to that of the DGFF (for the main idea supporting the proof we also refer to \cite{BDG}). We show here for the reader's convenience the upper bound, 
in order to obtain the desired limit in probability.\\
\noindent{\bf Proof of Theorem \ref{thm1.3} (a).} For any $\delta>0$ one can apply Chebyshev's inequality to get
\eqas
&&\prob\left(\left\{|\mathcal H_N(\eta)|\leq N^{-4(1-\eta^2)-\delta} \right\}\right)\leq  N^{4(1-\eta^2)+\delta} \av{|\mathcal H_N(\eta)|}\leq\\
&\leq & N^{-4(1-\eta^2)-\delta}N^4 \max_{x\in V_N}\prob\left(\varphi_x\geq2\sqrt{2 g}\eta \log N\right)\leq\\
&\leq & N^{-4(1-\eta^2)-\delta}N^4 \exp\left(-\frac{8 g\eta^2 \log^2 \! N}{2 g \log N+C}\right)\leq N^{-4(1-\eta^2)-\delta}N^{4-4\eta^2}\to 0
\eeqas
where we have used Lemma \ref{lemma2.1} too.
\endproof
\proofnam{thm1.4} We choose $\eta, \delta>0$ and define
\[
 D_+:=\left\{ \varphi_B\leq \cnst \eta \log N\right\},
\]
\[
 C_+:=\left\{|\hpts \cap D(x,N^\beta)|\geq N^{4 \beta(1-(\alpha/\beta)^2)-\delta)} \right\}
\]
and for an $\epsilon>0$ to be fixed later \[
 A:=\bigcup_{y \in B(x,N^\beta)}\left\{|\mathbf E(\varphi_y|\F_{\partial_2 B})-\varphi_B|\geq  \cnst\epsilon \log N \right\}.
\]
By Lemma \ref{lemma12} $\var(\varphi_B-\mathbf E(\varphi_y|\F_{\partial_2 B}))\leq c$ (we may assume that $B(x,N^\beta)\subsetneq V_N^{\ell}$), and so 
$$\prob(A)=O\left(N^{4 \beta}\exp\left(-c\log^2 \! N\right)\right)$$ 
tends to 0. Furthemore also $\prob(D^c_+)$ tends to $0$ by virtue of the bounds on covariances and \eqref{eq:2.6}. 
We then have
\begin{eqnarray*}
&& \prob(C_+)= \mathbf E(\prob (C_+|\F_{\partial_2 B}))\leq \prob(A)+\prob(D_+^c)+\mathbf E(\prob (C_+|\F_{\partial_2 B}){\mathbf{1}}_{A^c\cap D_+})\leq\\
&&\leq o(1)+\prob\left(\left|\mathcal H_{4N^\beta}\left(\frac{\alpha-\epsilon'}{\beta} \right)\right|\geq N^{4 \beta(1-(\alpha/\beta)^2)-\delta)}\right)
\end{eqnarray*}
where $\epsilon'$ satisfies
\[
 \frac{\alpha-\epsilon'}{\beta} \log( 4N^\beta)=(\alpha-\eta-\epsilon)\log N.
\]
By tuning the parameters $N$ large enough and $\eta, \epsilon$ small enough we can obtain
\[
 4\beta \left(1-\left(\frac{\alpha-\epsilon'}{\beta} \right)^2 \right)<4 \beta\left(1-\left(\frac{\alpha}{\beta} \right)^2 \right)+\delta
\]
(roughly speaking,  we have $\epsilon'\approx \alpha(1-\beta)$). By Theorem \ref{thm1.3} 
$$\prob\left(\left|\mathcal H_{4 N^\beta}\left(\frac{\alpha-\epsilon'}{\beta} \right)\right|\geq N^{4 \beta(1-(\alpha/\beta)^2)-\delta}\right)\to 0$$
and from this the claim follows. 
We now go to the lower bound proof, which is similar in spirit to the upper bound. By setting
\[
 D_-:=\left\{ \varphi_B\geq -\cnst \eta \log N\right\},
\]
\[
 C_-:=\left\{|\hpts \cap D(x,N^\beta)|\leq N^{4 \beta(1-(\alpha/\beta)^2)-\delta)} \right\}
\]
we also define
$$ \mathcal H_N^s(\eta):=\left\{ x\in V_N^s:\varphi_x\geq \cnst \eta\log N\right\},\quad s\in (0,1/2).$$
We observe that
\begin{eqnarray*}
&& \prob(C_-)= \mathbf E(\prob (C_-|\F_{\partial_2 B}))\leq \prob(A)+\prob(D_-^c)+\mathbf E(\prob (C_+|\F_{\partial_2 B}){\mathbf{1}}_{A^c\cap D_-})\leq\\
&&\leq o(1)+\prob\left(\left|\mathcal H^{3/8}_{4N^\beta}\left(\frac{\alpha+\epsilon'}{\beta}\right) \right|\leq N^{4 \beta(1-(\alpha/\beta)^2)-\delta)} \right)
\end{eqnarray*}
where $\epsilon'$ satisfies
\[
 \frac{\alpha+\epsilon'}{\beta} \log( 4N^\beta)=(\alpha+\eta+\epsilon)\log N
\]
and we conclude as before.\endproof
\proofnam{thm:cluster2} We will use the notation $b^\pm(\alpha,\beta,\eta,N)$ as in the proof of Lemma \ref{lemma8.1}. We will also introduce the following quantities: let $B:=B(x,4N^\beta)$, and for $\eta,\delta>0$,
\begin{eqnarray*}
E &:= & \left\{|\hpts \cap D(x,N^\beta)|\leq N^{4\beta(1-\alpha^2)-\delta} \right\},\\
F & := & \left\{\varphi_B\geq b^-(\alpha,\beta,\eta,N) \right\},\\
G &:=& \left\{x \in \hpts \right\}.
\end{eqnarray*}
\begin{description}
 \item[Lower bound.] Thanks to the proof of Lemma \ref{lemma8.1} we have $\prob(E|G)=\prob(E|F\cap G)\prob(F|G)+o(1)=\prob(E|F \cap G)(1+o(1))+o(1)$. This means that
\begin{eqnarray*}
&& \prob(E|F,G)=\frac{\prob(E\cap F \cap G)}{\prob(F\cap G)}\leq \frac{1}{\prob(F\cap G)}\sqrt{\prob(G)\prob(E\cap F)}=\\
&&= \frac{1}{\prob(F|G)\prob(G)}\sqrt{\prob(G)\prob (F)\prob(E|F)}=\\
&&\stackrel{Lemma\,\ref{lemma8.1}}{=}(1+o(1))\sqrt{\frac{\prob (F)}{\prob (G)}\prob(E|F)}.
\end{eqnarray*}
We know by the bounds \eqref{eq:var} and \eqref{eq:2.15}
\begin{eqnarray*}\prob(G)&=&\prob(\varphi_x \geq \cnst \alpha \log N)\geq c_1 \frac{\exp\left(-\frac{8g\alpha^2 \log^2 \! N}{2g \log N+c_2} \right)}{c_3\log N}\geq \\
& \geq &\exp\left(- d' \log N\right),\\
 \prob(F) & = & \prob(\varphi_B \geq \cnst (\alpha(1-\beta)-\eta) \log N)\leq \\
 &\leq & c_4 \exp\left(-\frac{8g(\alpha(1-\beta)-\eta)^2 \log^2 \! N}{2g(1-\beta) \log N+c_5} \right)\leq \exp\left(- d'' \log N\right)
\end{eqnarray*}
for some $d',\,d''>0$. Therefore we can find $d>0$ such that $\prob(F)/\prob(G)\leq \exp(d \log N)$ and to show the result it suf\/fices to prove that $\prob(E|F)\leq \exp(-c\log^2 \! N)$ for a positive $c$. For this purpose define
$$A:=\bigcup_{y \in B} \left\{\lvert\mathbf E(\varphi_y|\F_{\partial_2 B})-\varphi_B|\geq \cnst \epsilon \log N\right\}.$$
From Lemma \ref{lemma12} it follows that $\prob (A)\leq \exp(-c \log^2 \! N)$ 
for $c> 0$ and from \eqref{eq:2.15} that $\prob(F)\geq \exp\left(-d\log N\right)$ for some $d>0$, 
all in all $\prob(A|F)\leq \exp\left(-O\left( \log^2 \! N\right)\right)$. So we can write
\eqas &&
 \prob(E|F)\leq \frac{\prob(F\cap A)}{\prob(F)}+ \\
 &&+\frac{\prob(E\cap F \cap A^c)}{\prob(F)}\leq \\
 &&\exp\left(-O\left( \log^2 \! N\right)\right)+\frac{\mathbf E(\prob\left(E|\F_{\partial_2 B})\mathbf 1_{A^c}\mathbf 1_{F}\right)}{\prob(F)}.
\eeqas
If we are on $A^c \cap F$, then
\eqa &&
 \prob\left(|\hpts \cap D(x,N^\beta)|\leq N^{4 \beta (1-\alpha^2)-\delta}|\F_{\partial_2 B}\right)\leq \nonumber\\
 &&\leq\prob\left(\left|\mathcal H^{3/8}_{4N^\beta}(\alpha + \epsilon')\right|\leq N^{4 \beta (1-\alpha^2)-\delta}\right)\label{eq:4.1}
\eeqa
where $\epsilon'$ is such that 
\begin{equation}
 (\alpha- (\alpha(1-\beta)-\eta)+\epsilon)\log N=(\alpha+\epsilon')\log 4N^\beta.
\end{equation}
From Theorem \ref{thm1.3} we know that \eqref{eq:4.1} is bounded from above by $\exp(-c\log^2 \! N)$ for a constant $c>0$, provided that $\epsilon' $ is small 
(which can be obtained if 
$\eta$, $\epsilon$ and $N$ are small, small and large respectively).
\item[Upper bound.] Let $K \in \N$ and $\left\{\beta_j:=\frac{j}{K}\beta \right\}_{ 1\leq j\leq K}$. Then let
\[
 D_1:=D\left(x,N^{\beta_1}\right),\quad D_i:=D\left(x,N^{\beta_i}\right)\setminus D\left(x,N^{\beta_{i-1}}\right).
\]
Since $D\left(x,N^\beta\right)=\cup_{1\leq i\leq N}D_i$
\eqas &&
\left\{\left|\hpts \cap D\left(x,N^\beta\right)\right|\geq N^{4 \beta (1-\alpha^2)+\epsilon}  \right\}\subseteq \\
&& \subseteq\bigcup_{0\leq i\leq N}\left\{\left|\hpts \cap D_i\right|\geq N^{4 \beta_i (1-\alpha^2)+\epsilon/2} \right\}
\eeqas
as soon as $N$ is large. It is then suf\/ficient to prove that for all $i$
\[
 \prob\left(\left|\hpts \cap D_i\right|\geq N^{4 \beta_i (1-\alpha^2)+\epsilon/2} \left|x \in \hpts\right.\right)\stackrel{N\toi}{\longrightarrow}0.
\]
We can consider $\beta_j$'s for which $4 \beta_j(1-\alpha^2)+\epsilon/2\leq 4 \beta_j$. Let $B_j:=B\left(x,4N^{\beta_j}\right)$,
$$C:=\left\{\left|\hpts \cap D_j\right|\geq N^{4 \beta_j(1-\alpha^2)+\epsilon/2} \right\} $$ 
and $b^+(\alpha,\beta_j,\eta,N)$ as above. By Lemma \ref{lemma8.1} we obtain
\[\prob\left(C|x \in \hpts)=\prob(C \cap \left\{\varphi_{B_j}\leq b^+(\alpha,\beta_j,\eta,N)\right\}| x \in \hpts \right)+o(1).                                                                             
                                                                              \]
If we set $F:=\left\{\varphi_{B_j}\leq b^+(\alpha,\beta_j,\eta,N)\right\}$, $G:= \left\{x \in \hpts\right\} $ we obtain
\begin{eqnarray}
 &&\prob(C\cap F|G)\stackrel{Chebyshev\,inq.}{\leq}\frac{N^{-4\beta_j(1-\alpha^2)-\epsilon/2}}{\prob(G)}\mathbf E(\mathbf 1_{F\cap G}|\hpts \cap D_j|)=\nonumber\\
&&=\frac{N^{-4\beta_j(1-\alpha^2)-\epsilon/2}}{\prob(G)}\mathbf E \left(\sum_{y \in D_j}\indicator{x,y \in \hpts}\mathbf 1_{F} \right)\leq \nonumber\\
&&\leq \frac{N^{4\beta_j\alpha^2-\epsilon/2}}{\prob(G)}\sup_{y \in D_j}\prob(\left\{x,y \in \hpts\right\}\cap F).\label{eq:4.2}
\end{eqnarray}
By the bounds on the covariance and the normal distribution we have
\begin{equation}\label{eq:4.3}
 \prob(G)^{-1}\leq N^{4 \alpha^2+\epsilon/8}
\end{equation}
for $N$ large. By Lemma \ref{lemma8.3} by defining $\gamma^\ast=\frac{2}{2-\beta_j}>1$ when $\eta$ is small and $K$ large we obtain
%\begin{eqnarray}
 %&& 
\begin{equation}\sup_{y \in D_j}\prob(\left\{x,y \in \hpts \right\}\cap F)\leq  N^{-4\alpha^2 F_{2,\beta_j}(1)+\epsilon/8}
%&&
=N^{-4\alpha^2 (1+\beta_j)+\epsilon/8}.\label{eq:4.3bis}\end{equation}
%\end{eqnarray}
Inserting \eqref{eq:4.3} and \eqref{eq:4.3bis} in \eqref{eq:4.2} we obtain
\[\prob(C\cap F|G)\leq N^{4\beta_j \alpha^2-\epsilon/2+\epsilon/8 +4\alpha^2-4\alpha^2 (1+\beta_j)+\epsilon/8}=\frac{1}{N^{\epsilon/4}}\stackrel{N\toi}{\longrightarrow}0
\]                                                              
\end{description}
 \endproof 
 \proofnam{thm1.6} Preliminary we would like to make some considerations. It holds that $\rho(\alpha,\beta)$ is positive and in particular
\begin{equation} \rho(\alpha,\beta)\geq 4+4\beta -4\alpha^2 F_{2,\beta}(1)= 4(1-\alpha^2)(1+\beta)\label{eq:rho}.
\end{equation}
\eqref{eq:rho} derives from the fact that $F_{2,\beta}(\gamma)$ has a unique global minimum at $1$ in the range $\gamma\in\Gamma_{\alpha,\,\beta}$. Moreover notice that $\rho(\alpha,\beta)$ is increasing in $\beta$. If we set $\gamma_m:=\inf_{\gamma\in \Gamma_{\alpha,\beta}}F_{2,\beta}(\gamma)$, $\gamma_\ast:=\inf_{\gamma\geq 0}$ and $\gamma_+:=\sup{\Gamma_{\alpha,\beta}}$ we have $\gamma_\ast\leq \gamma_m\leq \gamma_+$ and moreover since $F_{h,\beta}(\cdot)$ does not depend on $\alpha$ as well as $\Gamma_{\alpha,\beta}$ does not depend on $\beta$ we have $\gamma_m=\min\left\{\gamma_\ast,\,\gamma_+\right\}$. 
 that
\[  \gamma_+=1/\alpha \geq 1 .\]
We are now ready to prove the lower and upper bounds.
\begin{description}
\item[Lower bound.]  We set
\[C:=\left\{\left|\left\{(x,y)\in \hpts:\,|x-y|\leq N^\beta \right\}\right|\leq N^{\rho(\alpha,\beta)-\delta}\right\}.
\]
Set $m_\gamma:=4-4\beta-4\alpha^2 F_{0,\beta}(\gamma)=4(1-\beta)(1-\alpha^2\gamma^2)$ and choose $\gamma<\gamma_+$ (in order to have $m_\gamma$ strictly positive). Further
\[F:=\left\{B \in \Pi_\beta:\,\varphi_B\geq \cnst \gamma(1-\beta)\alpha\log N \right\},
\]
\[D:=\left\{|F|\geq N^{m_\gamma-\delta/2} \right\}.
\]
Theorem \ref{thm1.3}\footnote{The idea is to scale the square: now we take the box with mesh $N/N^\beta$ and the grid is made by $\{x_{B}:\,B \in \Pi_\beta \}$. In this way Theorem \ref{thm1.3} tells us that $\mathcal H_{N^{1-\beta}}(\gamma \alpha)\approx N^{4(1-\beta)(1-\gamma^2 \alpha^2)}=N^{m_\gamma}$.} shows that $\prob(D^c)\to 0$. Hence we rewrite
\[\prob(C)=o(1)+\prob(D\cap C). \]
On $D$ we have at least $\left\{B_j:\,1\leq j\leq N^{m_\gamma-\delta/2} \right\}$ boxes. Set 
\[D_j:=\left\{\vr_{B_j}\geq \cnst \alpha \gamma(1-\beta)\log N \right\}.\]
We observe
\[
 C\cap D \subseteq E:=\bigcup_{j=1}^{N^{m_\gamma-\delta/2}} \left(D_j \cap \left\{\left|\hpts \cap B_j \right|\leq N^{(\rho(\alpha,\beta)-m_\gamma)/4-\delta/8}\right\} \right).
\]
Let us now put for some arbitrary $\eta>0$
\[
 A:=\bigcup_{B\in \Pi_\beta}\,\bigcup_{y\in B\left(x_B,N^{\beta}/2 \right)}\left\{ \left|\av\left(\varphi_y|\F_B\right)-\varphi_B\right|\geq \cnst \eta \log N\right\}.
\]
As before $\prob(A)=o(1)$ as $N\to+\infty$. Plugging this in, exactly as in the proof of Theorem \ref{thm1.4}
\eqas
 &&\prob(C\cap D)\leq o(1)+\prob(E\cap A^c)\leq\\
&&\leq o(1)+\\
&&+N^{m_\gamma-\delta/2}\prob\left(\left|\mathcal H_{N^\beta}^{1/4}\left(\frac{\alpha(1-\gamma(1-\beta))+\eta}{\beta}\right)\right|\leq N^{\frac{\rho(\alpha,\beta)-m_\gamma}{4}-\frac{\delta}{8}}\right).
\eeqas
Finally we observe that
\[
 \frac{\rho(\alpha,\beta)-m_\gamma}{4}\geq 2 \beta\left(1-\frac{\alpha^2\left(1-\gamma(1-\beta)\right)^2}{\beta^2}\right)
\]
which is $\exp(-O(\log^2 \! N))$ by Theorem \ref{thm1.3} for $\eta $ small enough, as we have already seen. Hence $\prob(C\cap D)=o(1)$, and we conclude the proof.
\item[Upper bound.] By Theorem \ref{thm1.3} we see that for $\lambda >0$ the number of $\alpha$-high points within distance $N^{\lambda \beta}$ is at most $N^{4(1-\alpha^2)+4\lambda \beta}$. We have with \eqref{eq:rho} that $4(1-\alpha^2)+4\lambda \beta\leq \rho(\alpha ,\beta)$ if
\[4(1-\alpha^2)+4\lambda \beta\leq 4(1-\alpha^2)(1+\beta)\iff \lambda
\leq (1-\alpha^2).\]
Therefore when this condition is not satisfied it is enough to find that there exists $h=h(\delta)<1$ such that for all $\beta'\in[\beta(1-\alpha^2),\beta]$
\[\prob\left(\left|\left\{(x,y)\in \hpts:\,N^{\beta'}\leq|x-y|\leq N^{\beta' h} \right\} \right|\geq N^{\rho(\alpha ,\beta')+\delta} \right) \rightarrow 0
.\]
We separate the two cases $\gamma_\ast\geqq \gamma_m$:
                           \begin{description}
                           \item[$\gamma_\ast= \gamma_m$.] Define 
                           \[E:= \left\{\left|(x,y)\in \hpts:\,N^{\beta'}\leq|x-y|\leq N^{\beta' h}  \right|\geq N^{\rho(\alpha ,\beta')+\delta}  \right\}.\]
                           By Chebyshev inequality
                           \eqas&&\prob(E)\leq N^{-\rho(\alpha ,\beta')-\delta} \mathbf E\left( \sum_{(x,y):\,N^{\beta'}\leq|x-y|\leq N^{\beta' h}}\right)\indicator{x,y\in \hpts} \leq\\
                           &&\leq N^{-\rho(\alpha ,\beta')-\delta} N^{4+4\beta'-4\alpha^2 F_{2,\beta'}(\gamma_\ast)+\delta/2},
                             \eeqas
                           where we have used the assumption that $h$ is close to 1 and Lemma \ref{lemma8.2}
                           \item[$\gamma_\ast> \gamma_m$.] We construct for each $B \in \Pi_{\beta'}$ a bigger box of size $4N^{\beta'}$ by juxtaposing to it the $12$ adjacent subboxes of same side length. We call the set of such bigger boxes $\mathcal B$, and for each $B' \in \mathcal B$ we center in $x_{B'}$ a box of twice bigger volume as $B'$. The latter boxes belong to a new set named $\mathcal C$. We remark that all pairs of points within distance $N^{\beta'}$ must belong to at least one $B' \in \mathcal B$. For $\epsilon>0$ set
                           \[ D:=\left\{\max_{C\in \mathcal C}\vr_C\geq (1+\alpha \epsilon)(1-\beta')\cnst \log N \right\}.
                            \]
                            By Lemma \ref{lemma2.1} and the fact that $\left\{\vr_y:\,y \in B\right\}$ with boundary conditions $\partial_2 B$ is a Gaussian field
                            \eqas \prob(D^c)&\leq& \left|\Pi_{\beta'} \right|\exp\left( -\frac{(1+\alpha\epsilon)^2(1-\beta')^2 (\cnst)^2\log^2 \! N}{2 g \log N^{\beta'}+O(1)}\right)\rightarrow 0
                            \eeqas
                            since $\left|\Pi_{\beta'} \right|=O(N^{4(1-\beta')})$. So noticing that $\alpha(\gamma_m+\epsilon)=(1+\alpha \epsilon)$
                            \[\prob(E)=o(1)+\prob(E \cap D)\leq o(1)+N^{-\rho(\alpha ,\beta')-\delta}N^{4+4\beta'-4\alpha^2 F_{2,\beta'}(\gamma_m+\epsilon)+\delta/2}
                            \]
                            if $h$ is close to $1$. $4+4\beta'-4\alpha^2 F_{2,\beta'}(\gamma_m+\epsilon)\stackrel{\epsilon\to 0}{\longrightarrow}\rho(\alpha ,\beta')$, thus $\prob(E)\to 0$.
                           \end{description}
\end{description}\endproof

\proofnam{thm:highsquare}
\begin{description}
\item[Lower bound.] We recall the notation used in the proof of Theorem \ref{thm1.3} by N. Kurt. For $\alpha \in (1/2,1)$ we choose $1\leq k \leq K+1$ such that
\eq \label{eq:6.1}\alpha_k:=\frac{\alpha(K-k+1)}{K}>\frac{1-\eta}{2}-\delta
\eeq
($\delta$ must be thought small). Let us now define recursively $\Gamma_{\alpha_1}:=\Pi_{\alpha_1}$. Then for $i\geq 2$, 
we set $\Gamma_{\alpha_i}$ as follows: for any $B\in \Gamma_{\alpha_{i-1}}$ define 
$\Gamma_{B,\alpha_i}:=\left\{B'\in \Pi_{\alpha_i} :\,B'\subseteq B/2\right\}$. Then 
$$\Gamma_{\alpha_i}:=\bigcup_{B\in \Gamma_{\alpha_{i-1}}}\Gamma_{B,\alpha_i}.$$ 
We re-use the notation $\underline B^{(k)}$ 
for a sequence of boxes $B_1 \supseteq B_2 \supseteq\cdots\supseteq B_k$, $B_i\in \Gamma_{\alpha_i}$ for all $1\leq i\leq k$. Finally
$$
D_k:=\left\{\underline B^{(k)}:\,\vr_{B_i}\geq (\alpha-\alpha_i)\lambda \cnst (1-1/K)\log N,\,\forall\,1\leq i\leq K \right\},
$$
$$
C_k:=\left\{\left|D_k\right|\geq n_k\right\}.
$$
We denote the biggest box of $\underline B^{(k)}$ with $B_{1,k}$. 
Let $B$ be a box of side length $N^{\alpha_k}/2$ centered in $B_{1,k}$. Let $n_k:=N^{\kappa+4\alpha(k-1)\frac{(1-\lambda)^2}{K}}$, where $\kappa$ is the 
constant appearing in \cite[Lemma 3.2]{Kurt_d5}. Define moreover for $\epsilon >0$
\[A:=\bigcup_{y \in B}\left\{ |\mathbf E(\varphi_y-\varphi_{x_B}|\F_{\alpha_k})|\geq \cnst \epsilon (\alpha-\alpha_k)(1-\gamma_K)\log N\right\}.
\]
By Lemma \ref{lemma12} $\prob(A^c)\rightarrow 1$ and $\prob(C_k)\rightarrow 1$ as in Theorem \ref{thm1.3} ($C_k$ is the same event). So
\eqas
&& \prob\left(D_N(\eta)\leq N^{\frac{1-\eta}{2}-\delta}\right)\leq\\
&& \leq o(1)+\\
&&+\prob\left(C_k \cap A^c\cap\left\{ \min_{y \in B}\varphi_y\geq\cnst \eta \log N\right\} \right)\leq\\
&& \stackrel{Def.\,of\,A,\,D_k}{\leq} o(1)+\\
&&+\prob\left( \min_{y \in B}(\varphi_y-\mathbf E(\varphi_y|\F_{\alpha_k}))\leq\right.\\
&&\left.\cnst \log N(\eta-(\alpha-\alpha_k)(1-\gamma_K)(1-\epsilon)) \right)\leq\\
&&\leq \prob\left( \max_{y \in V_{N^{\alpha_k}}^{1/4}}\varphi_y\geq\cnst \log N(-\eta+(\alpha-\alpha_k)(1-\gamma_K)(1-\epsilon)) \right)
\eeqas
where in the latter inequality we used the fact that $V_{N^{\alpha_k}}^{1/4}\supseteq B$.
For \eq \label{eq:sopra}\cnst \log N(-\eta+(\alpha-\alpha_k)(1-\gamma_K)(1-\epsilon))>\cnst \log N^{\alpha_k}\eeq
we would obtain thanks to Theorem \ref{thm1.3} that for $N$ large this probability tends to 0. But \eqref{eq:6.1} and \eqref{eq:sopra} give 
rise to a system of equations which has a solution for large $K$ and $N$, $\alpha$ close to 1 and $\epsilon$ small when $1/2+\eta/2k/K<\eta/2+\delta+1/2$.\endproof
\item[Upper bound.] We set $\theta:=\frac{1-\eta}{2}$, $\beta:=\theta+\delta$. We have first of all that
\eq \label{eq:6.2}\prob \left(\bigcup_{B \in \Pi_\beta}\{\varphi_B\geq \cnst (1-\theta )\log N\} \right)\stackrel{N \toi}{\longrightarrow}0
\eeq
since we have the variance bounds and \eqref{eq:2.6}. Furthermore let us define
\eqas
&& F:=\left\{\bigcap_{B \in \Pi_\beta} \{\varphi_B\leq \cnst (1-\theta )\log N\} \right\},\\
&& C:=\left\{ \bigcup_{B \in \Pi_\beta}\{\forall\,x \in B(\varphi_x\geq \cnst \eta\log N)\}\right\}.
\eeqas
We then have
\eqas  \prob\left(D_N(\eta)\right.&\geq &\left. N^{\theta+2 \delta}\right)\leq \prob(C)\leq\\
&&\leq \prob(F^c)+\prob(F \cap C)\leq \\
&& \stackrel{\eqref{eq:6.2}}{\leq} o(1)+\mathbf E(\prob(C|\F_\beta)\mathbf 1_F).
\eeqas
If $B \in \Pi_\beta$ we indicate with $B^{(1/4)}$ the sub-box $B(x_B,N^\beta/2)$. Choose $\epsilon >0$ and define
\[ A:=\bigcup_{B \in \Pi_B} \bigcup_{y \in B^{(1/4)}} \left\{|\mathbf E(\varphi_y -\varphi_{x_B}|\F_{\partial_2 B})|\geq \cnst \epsilon \log N \right\}.
\]
With Lemma \ref{lemma12} we obtain that $\prob(A)$ tends to $0$ as in Theorem \ref{thm1.4}. We can further bound
\[ \prob(D_N(\eta)\geq N^{\theta+2 \delta})\leq o(1)+\mathbf E(\prob(C|\F_\beta)\mathbf 1_{F\cap A^c}).
\]
To go on we notice that
\eq \label{eq:6.4}\prob(C|\F_\beta)\leq \left(\frac{N}{N^\beta} \right)^4 \max_{B \in \Pi_\beta}\prob\left(\forall\,x \in B(\varphi_x\geq \cnst \eta\log N)\right)
\eeq
and in particular on $F\cap A^c$
\eqas
&& \prob\left(\forall\,x \in B(\varphi_x\geq \cnst \eta\log N)\right) \leq\\
&& \leq  \prob\left(\forall\,x \in B(\varphi_x-\mathbf E(\varphi_x |\F_{\beta})\geq \cnst \log N(\eta-(1-\theta+\epsilon)))|\F_\beta \right)=\\
&&=\prob\left(\max_{x \in  V_{N^\beta}^{1/4}}\varphi_x \leq \cnst \log N(\theta+\epsilon)\right).
\eeqas
By Theorem \ref{Kurt1.2} this quantity is $O\left(\exp\left(-d \log^2 \! N\right)\right)$ for a positive $d$ when for instance $\beta>(\theta+\epsilon)$ which implies $\epsilon <\delta$. To sum up
\[ \prob(C|\F_\beta)\leq\exp\left(2(1-\beta)\log N-d \log^2 \! N\right)\rightarrow 0
\]
and recalling \eqref{eq:6.4} we finish the proof.\end{description}\endproof

\appendix
\section{Gaussian bounds}
\noindent{\bf Proof of \eqref{eq:2.6} and \eqref{eq:2.15}.\ }
\begin{description}
\item[$\eqref{eq:2.6}$] For $t>a>0$, $t+a>t-a$ and hence $t^2-a^2>(t-a)^2$,
\begin{eqnarray*}
\exp(a^2/2)P(|X|>a)&=&2\exp\left(a^2/2\right)P(X>a)=\\
&=&2\int_{a}^{+\infty}\frac{1}{\sqrt{2 \pi}}\exp\left(-\frac{t^2-a^2}{2}\right)\mathrm d t<\\
&< &2\int_{a}^{+\infty}\frac{1}{\sqrt{2 \pi}}\exp\left(-\frac{(t-a)^2}{2}\right)\mathrm d t=1.
\end{eqnarray*}
Notice that the bound holds also at a=0. 
\item[$\eqref{eq:2.15}$] We have that the function
$$
g(a):=2P(X>a)-\frac{\exp\left(-a^2/2\right)}{\sqrt{2 \pi}a}
$$
is such that $g(1)>0$, and its derivative
$$
g'(a)=\frac{2}{\sqrt{2 \pi}}\exp\left(-a^2/2\right)\left(\frac{1+a^2-a^3}{a^2}\right)<0,\quad \forall\,a\geq 1.
$$
Since $\lim_{a\to+\infty}g(a)=0$, $g(a)$ is always non negative.\end{description}

\section*{Acknowledgments} I would like to thank my supervisor Erwin Bolthausen for his guidance throughout this work. I am also grateful to the anonymous referee for his careful review of an earlier draft and for suggestions for improvement.

\bibliography{literatur}
\end{document}